\documentclass[11pt]{article}
\usepackage{latexsym,amssymb}
\def\sym{\mathbb}
\def\F{{\sym F}}
\def\K{{\sym K}}
\def\N{{\sym N}}
\def\Z{{\sym Z}}
\def\Q{{\sym Q}}
\newtheorem{prop}{Proposition}
\begin{document}
\begin{center}
{\LARGE\bf Extending the scalars of minimizations
}\\[5mm]
{\large G. DUCHAMP\footnote{Gerard.Duchamp@univ-rouen.fr}\hspace{1cm}
\'E. LAUGEROTTE\footnote{Eric.Laugerotte@univ-rouen.fr}
}\\[1mm]
{\bf L}aboratoire d'{\bf I}nformatique {\bf F}ondamentale et 
{\bf A}ppliqu\'ee de {\bf R}ouen\\
Facult\'e des Sciences et des Techniques\\
76821 Mont Saint Aignan CEDEX France\\[2mm]
{\large J-G. LUQUE\footnote{Jean-Gabriel.Luque@univ-mlv.fr\\
Partially supported by the Scientific Research Program of MENRT}}
\\[1mm]
{\bf I}nstitut {\bf G}aspard {\bf M}onge\\
Universit\'e de Marne la Vall\'ee\\
77454 Marne la Vall\'ee France
\end{center}

\section{Introduction}
In the classical theory of formal languages, finite state automata allow to 
recognize the words of a rational subset of
$\Sigma^*$ where $\Sigma$ is a set of symbols (or the alphabet). Now, 
given a semiring $(\K,+,.)$, one 
can construct $\K$-subsets of $\Sigma^*$ in the sense of Eilenberg \cite{Ei}, 
that are alternatively
called noncommutative formal power series \cite{BR,St} for which a framework
very similar to language theory has been constructed (see \cite{Sc1,Sc2} and
\cite{BR}). This extension has applications in many domains. 
Let us cite, for example, enumeration 
(non-commutative as used by instance for alignment of genomic sequences), 
image processing \cite{CK}, task-ressource problems \cite{GM} 
and real-time applications where multiplicities are used to
prove the modularity of the validation method by means of the Hadamard product
of two integer valued automata (see the contribution by Geniet and Dubernard 
\cite{GD}).

Particular noncommutative formal power series, which are called rational 
series, are the behaviour of a family of weighted automata 
(or $\K$-automata). In order to get an efficient encoding, 
it may be interesting to point out one of them with the smallest number of 
states. 
Minimization processes of $\K$-automata already exist for $\K$ being:\\ 
{\bf a)} a field \cite{BR},\\
{\bf b)} a noncommutative field \cite{FL},\\
{\bf c)} a PID \cite{Fl}.\\
When $\K$ is the bolean semiring, such a minimization process (with 
isomorphisms of minimal objects) is known within the category of deterministic 
automata. 

Minimal automata have been proved to be isomorphic in cases {\bf (a)} and 
{\bf (b)} (see respectively \cite{BR} and \cite{FL}). 
The case {\bf (c)} is mentioned in \cite{Fl}. But the 
proof given in \cite{BR} is not constructive. In fact, it lays on the 
existence of a basis for a submodule of $\K^n$. Here we give an independent 
algorithm which reproves this fact and an example of a
pair of nonisomorphic minimal automata. Moreover, we examine the possibility 
of extending {\bf (c)}. To this end, we provide
an {\em Effective Minimization Process} (or {\em EMP}) which can be used for 
more general sets of coefficients. 

The structure of the contribution is the following. After this introduction,
we give in
details the EMP and, in Section 3, we discuss the termination of the EMP in a 
frame which extends {\bf (c)}. 

\section{Computing a prefix subset}
Let $\K$ be an integral domain (a ring without zero divisor) and 
$\Sigma$ a finite alphabet. 
A $\K$-automaton $\cal A$ is usually identified by a linear representation 
$(\lambda,\mu,\gamma)$. We examine here a process which allows us to find a 
prefix subset $X$ such that $\lambda\mu(X)$ generates 
$\lambda\mu(\K\langle A\rangle)$. We apply Algorithm {\bf prefix} (which calls 
Algorithm {\bf generator}) to a $\K$-automaton $\cal A$. 

\bigskip
\noindent
{\bf Algorithm prefix}\\
\begin{tabular}{lll}
input&:&the linear representation $(\lambda,\mu,\gamma)$.\\
output&:&a pair $(X,Z)$ where $X$ is a prefix code and $Z\subset X$.\\
\end{tabular}
\begin{enumerate}
\item $(X_0,Y_0,Z_0):=(\emptyset,\{1\},\emptyset)$
\item if $Y_i\neq\emptyset$
\begin{enumerate}
\item choose $y\in Y_i$ of minimal length
\item $(X_{i+1},Y_{i+1},Z_{i+1}):=\mbox{\bf
generator}((\lambda,\mu,\gamma),y,
(X_i,Y_i,Z_i))$
\item go to (2)
\end{enumerate}
\item return $(X,Z)$
\end{enumerate}

\bigskip
\noindent
{\bf Algorithm generator}\\
\begin{tabular}{lll}
input&:&the linear representation $(\lambda,\mu,\gamma)$,\\
&&the word $y$,\\
&&the triplet $(X,Y,Z)$.\\
output&:&the triplet $(X,Y,Z)$.\\
\end{tabular}
\begin{enumerate}
\item $n:=|X|$
\item if it does not exist $\alpha\in\K$ such that
$\alpha\lambda\mu(y)=\alpha_1\lambda\mu(x_1)+\cdots
+\alpha_n\lambda\mu(x_n)$
with $\alpha_i\in K$ and $x_i\in X$ ($1\leq i\leq n$) then
\begin{itemize}
\item[] $(X,Y,Z):=(X\cup\{y\},Y\cup yA\setminus \{y\},Z)$
\end{itemize}
\item else if it exists such a $\alpha$ which divides $\alpha_i$ 
($1\leq i\leq n$) 
\begin{itemize}
\item[] $(X,Y,Z):=(X,Y\setminus \{y\},Z)$
\end{itemize}
\item else
\begin{itemize}
\item[] $(X,Y,Z):=(X\cup\{y\},Y\cup yA\setminus \{y\},Z\cup \{y\})$
\end{itemize}
\item return $(X,Y,Z)$
\end{enumerate}

\bigskip
\noindent
As a computation process \cite{Kn}, Algorithm {\bf prefix} is well-defined if 
we can 
compute $\alpha$ and the $\alpha_i$'s in Algorithm {\bf generator}. 
Let us denote $\F$ the field of fractions of $\K$. 

\begin{prop}
When the computation process terminates, 
\begin{enumerate}
\item the family $\lambda\mu (X)$ generates $\lambda\mu(\K\langle
A\rangle )$.
\item the family $\lambda\mu (X-Z)$ is linearly independent for $\F$.
\end{enumerate}
\end{prop}
We prove 1 as in \cite{BR} using the decomposition $\Sigma^*=C^*X$ with 
$C=(X\Sigma\cup\epsilon)\setminus X$ the prefix code induced by $X$, 
and using the linearity of $\mu$. Now, for 2, the only way to make the set 
$X-Z$ increasing is to come through Step 2 of Algorithm {\bf generator} where 
we add to $X$ an item $y$ such
that $\lambda\mu(y)$ is linearly independant of $\lambda\mu(X-Z)$.

\bigskip
The family $\lambda\mu (X-Z)$ does not generate necessarily
$\lambda\mu(\K\langle A\rangle)$ but we could expect that it exists a basis of
$\lambda\mu(\K\langle A\rangle )$ of rank $|X-Z|$. This occurs only when the 
relation $\alpha \lambda\mu(y)=\sum_{x\in X}\alpha_x\lambda\mu(x)$ 
implies that the rank of $\mbox{Span}(\lambda\mu(X\cup\{y\}))$ is $|X|$. Or 
again, this is
equivalent to respect the following condition: 
for each $n,\ m\in \N^+$, if $V=\{v_i\}_{i\in[1,m]}\subseteq \K^n$ is
linearly independent and $\alpha u=\sum_i \alpha_i v_i$ (with $\alpha
\in\K-\{0\}$) then the rank of $\mbox{Span}(V\cup\{u\})$ is $m$. Then, 
setting $m=1$ and $n=1$, we find that such a ring $\K$ verifies the B\'ezout 
condition. We will say that $\K$ is a B\'ezout ring. 
Conversely, suppose that $\K$ is a B\'ezout ring. Using a Gauss
method, we find the property. More precisely, let $\left(a\atop
b\right)\in \K^{2}$. If $b=0$ the triangularization is clear. If $a=0$ then
\[\left(
\begin{array}{cc}
0&1\\1&0
\end{array}\right)\left(0\atop b\right)=\left(b\atop0\right).\]
Otherwise, as $\K$ is a B\'ezout ring, it exists $\alpha,\beta\in\K$ such
that
\[\alpha a+\beta b=\gcd(a,b)=d.\]
Then, if the matrix $G$ is defined by \[G=\left(\begin{array}{cc}
\alpha&\beta\\
-\displaystyle\frac bd&\displaystyle\frac ad\end{array}
\right) \]
one has
\[ G\left(a\atop b\right)=\left(d\atop 0\right).\]
But, as $\K$ is an integral domain, the matrix $G$ is unimodular. 
We can apply this process to triangularize matrices in $\K^{n\times n}$. 
In the sequel, we will only consider an integral B\'ezout
domain $\K$.

\section{Minimization in integral B\'ezout domains}
We use the prefix code computed in Algorithm {\bf prefix} in order to
construct a left reduced of a linear representation $(\lambda,\mu,\gamma)$.
The main step of our algorithm is
to choose a basis of $\lambda\mu\left(\K\langle A\rangle\right)$ using the
previous Gauss process. In fact, we consider Algorithm {\bf triang} 
taking a matrix $M$ as input and returning a pair $(G,T)$ where $T$ is a
stair matrix and $G$ a Gauss matrix such that 
\[
\left(\frac
T{0^n}\right)=\left(\begin{array}{cc}G&0\\0&Id_n\end{array} \right)M
\]
with $n$ maximal.

\eject
\noindent
{\bf Algorithm left\_reduction}\\
\begin{tabular}{lll}
input&:&a linear representation $(\lambda,\mu,\gamma)$.\\
output&:&a left reduced linear representation
$(\lambda_r,\mu_r,\gamma_r)$.\\
\end{tabular}
\begin{enumerate}
\item $(X,Z):=\mbox{\bf prefix}((\lambda,\mu,\gamma))$
\item $(I,T)=((1),(\lambda))$
\item if $X\not =\emptyset$ then
\begin{enumerate}
\item choose $x\in X$ of minimal length
\item $X:=X\setminus \{x\}$
\item $(G,T):=\mbox{\bf triang}\left(\left(\frac
T{\lambda\mu(x)}\right)\right)$
\item if $x\in Z$ then $I:=IG^{-1}$ else $I:=(I|0)G^{-1}$
\item go to (3)
\end{enumerate}
\item\label{item2} for each $a\in A$, compute $\mu_r(a)$ such that
$T\mu(a)=\mu_r(a)T$.
\item return $(I,\mu_r,T\gamma)$
\end{enumerate}

\bigskip
\noindent
In Step 3(d), the ``if'' part occurs when the rank of the matrix 
$\left(\left(\frac T{\lambda\mu(x)}\right)\right)$ with coefficients in $\Z$ 
is equal to the rank of $T$. In the ``else'' part, we add an item to the 
family, and then one line to $T$ and one column to I. 

\begin{prop}
Let ${\cal A}=(\lambda,\mu,\gamma)$ be a linear representation. 
When the computational method terminates, 
Algorithm {\bf left\_reduction} gives a left reduced $\K$-automaton of ${\cal A}$. 
\end{prop}

\noindent
We can observe that 
$\lambda\mu(w)\gamma=\lambda_rT\mu(w)\gamma=\lambda_r\mu_r(w)T\gamma=\lambda
_r\mu_r(w)\gamma_r$. 
Furthermore, the construction implies that the linear representation 
$\lambda_r\mu_r(\F\langle X\rangle )$ lies in $\F^{1\times |X-Z|}$. 
Moreover, we can compute a right reduced automaton using the previous 
algorithm with the linear
representation $(\gamma^t,\mu^t,\lambda^t)$ as input. Realizing a left 
reduction and a right reduction gives a minimal $\K$-automaton. 
However, here, two minimal linear representations are not necessarily 
isomorphic. As shown 
the following example: 
\begin{enumerate}
\item ${\cal A}_1=\left((1\quad
0),\left(\begin{array}{cc}0&x\\0&0\end{array}\right),
\left(\begin{array}{c}0\\1\end{array}\right)
\right)$
\item ${\cal A}_2=\left((x\quad
0),\left(\begin{array}{cc}0&1\\0&0\end{array}\right),
\left(\begin{array}{c}0\\1\end{array}\right)
\right)$
\end{enumerate}
with $x\in \Z-\{0,1,-1\}$. These $\K$-automata are different minimal 
linear representations with a same behaviour, but they are not isomorphic. 
Let $T$ be a matrix such that
\[{\cal A}_1=\left((x\quad 0)T^{-1}
,T\left(\begin{array}{cc}0&1\\0&0\end{array}\right)
T^{-1},T\left(\begin{array}{c}0\\1\end{array}\right)\right). 
\]
This relation implies necessarily that 
\[(T^{-1})_{2,1}=\frac1x\not\in\Z.\]

\bigskip
\noindent
{\bf Note}\hspace{2mm} 
Let $\F$ be a field, then 
$\K=\F((X^\alpha)_{\alpha\in \Q_+\setminus\{0\}})$ 
(polynomials with fractional powers) provides an example of integral B\'ezout 
domain where Algorithm {\bf left\_reduction} terminates, and which is not a 
principal integral domain.

\end{document}